\pdfoutput=1
\RequirePackage{ifpdf}
\ifpdf 
\documentclass[pdftex]{sigma}
\else
\documentclass{sigma}
\fi

\begin{document}

\allowdisplaybreaks

\renewcommand{\thefootnote}{$\star$}

\renewcommand{\PaperNumber}{063}

\FirstPageHeading

\ShortArticleName{Leibniz Algebras and Lie Algebras}

\ArticleName{Leibniz Algebras and Lie Algebras\footnote{This paper is a~contribution to the Special Issue on New Directions in Lie Theory.
The full collection is available at \href{http://www.emis.de/journals/SIGMA/LieTheory2014.html}
{http://www.emis.de/journals/SIGMA/LieTheory2014.html}}}

\Author{Geoffrey MASON~$^\dag$ and Gaywalee YAMSKULNA~$^{\ddag\S}$}

\AuthorNameForHeading{G.~Mason and G.~Yamskulna}

\Address{$^\dag$~Department of Mathematics, University of California, Santa Cruz, CA 95064, USA}
\EmailD{\href{mailto:gem@ucsc.edu}{gem@ucsc.edu}}

\Address{$^\ddag$~Department of Mathematical Sciences, Illinois State University, Normal, IL 61790, USA}
\EmailD{\href{mailto:gyamsku@ilstu.edu}{gyamsku@ilstu.edu}}

\Address{$^\S$~Institute of Science, Walailak University, Nakon Si Thammarat, Thailand}

\ArticleDates{Received September 09, 2013, in f\/inal form October 19, 2013; Published online October 23, 2013}

\Abstract{This paper concerns the algebraic structure of
f\/inite-dimensional complex Leibniz algebras. In
particular, we introduce left central and symmetric
Leibniz algebras, and study the poset of Lie subalgebras
using an associative bilinear pairing taking values
in the Leibniz kernel.}

\Keywords{Leibniz algebras; Lie algebras}

\Classification{17A32}

\renewcommand{\thefootnote}{\arabic{footnote}}
\setcounter{footnote}{0}

\section{Introduction}

Throughout the present paper, a \emph{left Leibniz algebra} means a nonassociative $\mathbb{C}$-algebra $M$ with a  product  (or \emph{bracket}) $[ \ ]$ satisfying the following identity for all $a, b, c \in M$ :
\begin{gather}
[a[bc]]=[[ab]c]+[b[ac]]. \label{rderivid}
\end{gather}
 (\ref{rderivid}) says that the left adjoint ${\rm ad}_a: b \mapsto [ab]$ $(a, b \in M)$ is a \emph{derivation} of  $M$,
so that  ${\rm ad}:   M \rightarrow {\rm  Der}(M)$, $a \mapsto {\rm ad}_a,$ is a morphism
of $M$ to the algebra of derivations of $M$ regarded as a Lie (or left Leibniz) algebra with the usual bracket. The \emph{Leibniz kernel} is the subspace $C(M)$ spanned by $[aa]$ $(a \in M)$; $M$ is a Lie algebra if, and only if, $C(M)=0$.

Leibniz algebras were f\/irst studied for their own sake by Loday \cite{L1} (see also \cite[Section~10.6]{L2}). The rationale for the present
work is partially motivated by the  triangular decomposition
\begin{gather}\label{triangdecomp}
V =  (\oplus_{n \leq 0} V_n  )\oplus V_1 \oplus  (\oplus_{n \geq 2} V_n )
\end{gather}
of a vertex operator algebra (VOA) $V$. In the most widely studied case when~$V$ is of
\emph{CFT-type}, i.e.\ $V_n=0$ for $n<0$ and $V_0=\mathbb{C}\mathbf{1}$ is spanned by the vacuum vector, the summands in~(\ref{triangdecomp}) are Lie algebras.  In the general case~(\ref{triangdecomp}) satisf\/ies only the weaker condition of being a~decomposition into $\mathbb{Z}$-graded \emph{left Leibniz algebras} with respect to
the $0^{\rm th}$ operation in $V$ (cf.~\cite{LL}). This decomposition plays a r\^{o}le in our work~\cite{MY}
where, among other things, we are interested in the Lie subalgebras of $V_1$ and the vertex subalgebras
of $V$ that they generate.  This leads directly to the main theme of the present paper, which is the study of the poset of Lie subalgebras of certain kinds of Leibniz algebras.  Readers uninured to VOA theory need not be concerned~--  the present paper deals solely with Leibniz algebras, and no further mention of VOAs will be made.

The def\/inition of Leibniz algebra is \emph{not} left-right symmetric; a right Leibniz algebra, in which
the map $b \mapsto [ba]$ for f\/ixed $a$ is a derivation, is not necessarily a left Leibniz algebra. We call an algebra that is both a left and right Leibniz algebra a \emph{symmetric} Leibniz algebra.
Notice from~(\ref{rderivid}) that a left Leibniz algebra $M$ satisf\/ies
 \begin{gather}\label{nearcommid}
[[aa]b]=0, \qquad a, b \in M,
\end{gather}
and dually, a right Leibniz algebra satisf\/ies
\begin{gather}\label{centralid}
[a[bb]]=0, \qquad a, b \in M.
\end{gather}
We call~$M$ a left \emph{central} Leibniz algebra if it is a left  Leibniz algebra that also satisf\/ies~(\ref{centralid}). Equivalently, $M$ is both the left and right centralizer of~$C(M)$.
There is a hierarchy of algebras
\begin{gather*}
 \{\mbox{left Leibniz}\} \supseteq \{\mbox{left central Leibniz}\}  \supseteq \{\mbox{symmetric Leibniz}\}
\supseteq \{\mbox{Lie}\},
\end{gather*}
and in fact each containment is \emph{strict}.

We now describe the contents of the present paper, which deals solely with
\emph{finite-dimensional}, complex, left Leibniz algebras~$M$. We tacitly assume this in everything that follows.
In Section~\ref{section2} we discuss some basic facts, in particular concerning Levi subalgebras and
Levi decompositions of Leibniz algebras. Levi subalgebras (i.e., semisimple Lie subalgebras that complement the solvable radical) are readily seen to exist in~$M$ (see~\cite{B1} and Section~\ref{section2} below). Malcev's theorem for Lie algebras
does \emph{not} extend to left Leibniz algebras, though it does for various special classes,
including left central Leibniz algebras. Section~\ref{section3} is devoted to these issues.

The remainder of the paper revolves about the symmetric bilinear pairing $\psi: M\times M \rightarrow C(M)$,
$(a, b) \mapsto [ab]+[ba]$. This is a feature of all Leibniz algebras, but for left central Leibniz algebras
it is particularly ef\/f\/icacious, because in this case it is  \emph{associative}. Then the radical $R$ of $\psi$ is a $2$-sided ideal of~$M$,
and the poset $\mathcal{L}$ of Lie subalgebras of~$M$ coincides with the set of Leibniz subalgebras which are also isotropic subspaces. General properties of this set-up are developed in Section~\ref{section4}, including (Proposition~\ref{prop4.1}) the fact that $R$ coincides with the intersection of the
maximal elements of $\mathcal{L}$. We also give (Lemma~\ref{lemtildeM}) a general construction of a class of left central Leibniz algebras based solely on Lie-theoretic data. In Sections~\ref{section5} and~\ref{section6} we consider left central Leibniz algebras with $\dim C(M)=1$ (so that $\psi$ is a trace form in the usual sense). Such algebras arise, for example, as quotients~$M/U$ whenever~$U$ is a hyperplane of~$C(M)$. We prove (Theorem~\ref{thm5.1}) that
there is a  Lie subalgebra $L\subseteq M$ that is also a maximal isotropic subspace of~$M$
(and therefore also a maximal element of~$\mathcal{L}$) and such that~$L/R$ is \emph{nilpotent}.
In Section~\ref{section6} we assume that $M$ is also symmetric, a property that holds precisely when $M'\subseteq R$ (Lemma~\ref{lemLsymm}).
 We prove (Theorem~\ref{thm6.3}) that a symmetric Leibniz algebra with $\dim C(M)=1$ has a $2$-sided ideal of codimension at most 1 that arises from the construction of Lemma~\ref{lemtildeM}. In this way, we more-or-less obtain a characterization of such Leibniz algebras  in terms of Lie-theoretic data.

\section{Basic facts}\label{section2}

In this section, $M$ is a left Leibniz algebra.
For subsets $A, B \subseteq M$, we def\/ine $[AB]$ (or $[A, B]$) to be the subspace
spanned by all brackets $[a b]$ $(a\in A$, $b\in B)$. If $A=\{a\}$ is a singleton, we write
$[\{a\} B] = [a, B]$.
Introduce
\begin{gather*}
 Z(M) := \{z \in M \, | \, [M, z] = [z, M] = 0\}, \\
  M' := [M M], \\
  C(M) := \mbox{span} \langle [aa]  \, | \, a \in M \rangle = \mbox{span}\langle [ab]+[ba] \, | \, a, b \in M \rangle.
\end{gather*}
$Z(M)$, $M'$ and $C(M)$ are $2$-sided ideals of $M$ called the \emph{center}, \emph{derived subalgebra} and
\emph{Leibniz kernel} of~$M$ respectively. (The second equality def\/ining~$C(M)$ is equivalent to the f\/irst by polarization. That $C(M)$ is a $2$-sided ideal then follows from (\ref{nearcommid}) and (\ref{rderivid}) with $b=c$. See also~\cite{AO,B2}.)  Obviously,
$M$ is a Lie algebra if, and only if, $C(M)=0$.

Set $M^{(0)}:= M$ and for $n \geq 0$ inductively def\/ine $M^{(n+1)}:= (M^{(n)})'$;
$M$ is called \emph{solvable} if $M^{(n)}=0$ for some $n\geq 0$. $M$~has a unique maximal solvable
ideal, called the \emph{solvable radical} of~$M$ and denoted by~$B(M)$. $M$~is called \emph{abelian} if~$M'=0$.
By (\ref{nearcommid}), $M$ is the right centralizer of~$C(M)$, i.e.,
\begin{gather}\label{nearLA}
[C(M)M]=0.
\end{gather}
Consequently,  $C(M)$
is a $2$-sided abelian ideal of $M$ and $C(M)\subseteq B(M)$. Because $C(M/C(M))$ $=0$, $M/C(M)$ is a Lie algebra and $B(M)/C(M)$ is its solvable radical. $M$~is a left central Leibniz algebra if, and only if,
$C(M)\subseteq Z(M)$.

It is not hard to see that the analog of the \emph{Levi decomposition} for Lie algebras holds in $M$
(cf.~\cite{B1}). Indeed, set $N=C(M)$, $B=B(M)$, and let $N\subseteq T \subseteq M$  with $T/N$ a Levi subalgebra of the Lie algebra~$M/N$. We have a Levi decomposition in~$M/N$,
\begin{gather*}
M/N = B/N \oplus T/N.
\end{gather*}
 Because $T' \subseteq T$ and $N$ acts trivially on the left of~$M$, $T$ is a left $T/N$-module. Because $T/N$ is semisimple, Weyl's theorem of complete reducibility tells us that there is a left $T$-submodule that complement~$N$ in~$T$, call it $S$. Therefore, $S' \subseteq [T S] \subseteq S$, so that $S$ is a Leibniz subalgebra of~$M$. Moreover,
$C(S) \subseteq S \cap N=0$, whence~$S$ is, in fact, a Lie subalgebra of~$M$. It is clearly isomorphic to~$T/N$, hence is semisimple. We have a direct sum decomposition
\begin{gather}\label{Levidecomp}
M = B \oplus S.
\end{gather}
We call a Lie subalgebra $S$ of~$M$ that complements~$B$ as in (\ref{Levidecomp}) a \emph{Levi factor}, or \emph{Levi subalgebra}, of~$M$;  (\ref{Levidecomp}) itself is called a \emph{Levi decomposition} of~$M$. The Levi
subalgebras of~$M$  can be characterized as the Lie subalgebras of~$M$ of maximal dimension subject to being semisimple.

\section{On Malcev's theorem for Leibniz algebras}\label{section3}

Malcev's theorem for f\/inite-dimensional complex Lie algebras includes the statements that all Levi subalgebras are conjugate by the exponential of an inner derivation, and
every semisimple Lie subalgebra is contained in a Levi subalgebra. We shall see that both of these assertions are generally false for Leibniz algebras. On the other hand, versions of the Malcev theorem can be proved for certain classes of Leibniz algebras. This section is concerned with these questions.

We begin with a general construction. Let~$S$ be a f\/inite-dimensional complex Lie algebra and~$N$ a
f\/inite-dimensional left $S$-module with action $S\times N \rightarrow N, (a, m) \mapsto a.m$ $(a \in S$, $m \in N)$. Let $M=S \oplus N$ with bracket
\begin{gather}\label{firstbrack}
[(a, m)(b, n)] := ([ab], a.n).
\end{gather}
Here, of course, $[ab]$ denotes the bracket in~$S$. (\ref{firstbrack}) def\/ines the structure of
a Leibniz algebra on~$M$; $S$~is a Lie subalgebra of~$M$,
and $C(M)= [SN]$.

Now suppose in addition that $S$ is semisimple. Then $S$ is a Levi subalgebra
of~$M$ and~$N$ is the solvable radical. Suppose that $S_1$ is any Levi subalgebra of~$M$.
Then $M=S_1\oplus N$ and $S_1 =\{ (x, c(x))\, | \, x \in S\}$ for some linear map
$c: S \rightarrow N$. Moreover
$[(x, c(x))(y, c(y))] = ([xy], x. c(y)) = ([xy], c([xy])),$ so that
$x. c(y) = c([xy])$. This says that
$c$ is a morphism of left $S$-modules, where~$S$
furnishes the left adjoint $S$-module. Conversely, given
$c \in {\rm Hom}_S(S, N)$, the set of pairs $\{(x, c(x)) \, | \, x \in S\}$ is a Levi subalgebra of~$M$.

Suppose that $S$ is  simple and
$N$ a simple $S$-module that is \emph{not} the adjoint module. Then ${\rm Hom}_S(S, N)=0$,
so that $S$ is the \emph{unique} Levi subalgebra of~$M$. In the case of Lie algebras, if there is a unique Levi subalgebra then it is an ideal. However, as our construction shows, this is not
necessarily true for Leibniz algebras.

Suppose that $T\subseteq S$ is a simple Lie subalgebra. Replacing $S$ by $T$ in the previous paragraphs, we see that if
$d \in {\rm Hom}_T(T, N)$ then $T_1 := \{(x, d(x)\, | \, x \in T\}$ is a Lie subalgebra of $M$ isomorphic
to $T$. In order for $T_1$ to be contained in a Levi subalgebra, it is necessary and suf\/f\/icient
that $d$ is the \emph{restriction} $d = {\rm Res}^S_{T}(c)$ for some $c \in {\rm Hom}_S(S, N)$.

\begin{example}\label{example3.1}
 $S= \frak{sl}_3$, $V$ the natural $3$-dimensional module,
$N=S^2(V)$ the second symmetric square, $T$ the $\frak{sl}_2$-subalgebra of $S$ corresponding to a simple root. Then ${\rm Hom}_S(S, N)=0$ and ${\rm Hom}_T(T, N)\not= 0$. Thus there exists $T_1$ \emph{not} contained in the unique Levi subalgebra~$S$.
The theory of highest weight modules can be used to construct many similar examples.
\end{example}

\begin{example}\label{example3.2}
 Suppose that $S$ is simple and there is an isomorphism of left $S$-modules
$S \cong N$. Then ${\rm Hom}_S(S, N)\not = 0$, and there are at least two distinct Levi subalgebras
in~$M$. On the other hand, because $[NM]=0$, every inner derivation of~$M$ leaves every
Levi subalgebra invariant, so that conjugacy \emph{a la} Malcev does not hold in~$M$. (This is
essentially Example~2 in~\cite{B1}.)
\end{example}

Now we turn to some positive results in the direction of Malcev's theorem.

 \begin{theorem}\label{lemmaconj1}
 Let $M$ be a left Leibniz algebra
 with Leibniz kernel $N$ and solvable radical~$B$. Let~$S$ be a Levi subalgebra of~$M$, and suppose that $[S N]=0$. If $T\subseteq M$ is any semisimple Lie subalgebra, there
 is an  automorphism $\alpha=\exp({\rm ad}\, x)$ $(x \in [M, B])$
 such that $\alpha(T)\subseteq S$.
  \end{theorem}

\begin{remark}
 If one Levi subalgebra $S$ satisf\/ies $[SN]=0$, then
 \emph{all} Levi subalgebras have the same property.
 \end{remark}

\begin{proof}[Proof of Theorem~\ref{lemmaconj1}.] Let $T\subseteq M$ be any semisimple
Lie subalgebra. By Malcev's theorem for Lie algebras applied to~$M/N$, we can f\/ind
$x+N \in [M/N, B/N]$ so that
$\alpha =\exp({\rm ad}\, (x+N))$ and $\alpha(T+N/N) \subseteq S+N/N$. Because
of~(\ref{nearLA}), we can identify~$\alpha$ with
$\exp({\rm ad}\, x)$, which we also denote by $\alpha$. Then $\alpha(T)\subseteq S+N$.
Because $[M/N, B/N] = ([M, B]+N)/N$, we may choose $x \in [M, B]$.
Because of the assumption that $[S N]=0$,
$S+N$ is a \emph{reductive} Lie subalgebra, and~$S$ is the unique
maximal semisimple subalgebra. Therefore $\alpha(T)\subseteq S$, and
the theorem is proved.
\end{proof}

\begin{theorem}\label{lemmaconj3} Let $M$ be a left Leibniz algebra
 with Leibniz kernel~$N$ and solvable radical~$B$. Suppose that $[BN]=0$.
 Then given any pair of Levi subalgebras~$S_1$ and~$S_2$  of~$M$, there is a~derivation $f$ of $M$ such that $\exp(f)(S_1)=S_2$.
\end{theorem}

\begin{remark}
 The conditions of the theorem apply in the situation of Example~\ref{example3.2}, for example.
Thus, while we may not have conjugacy of Levi subalgebras using only exponentials of inner derivations, conjugacy is restored in some cases if more general automorphisms are permitted.
\end{remark}

\begin{proof}[Proof of Theorem~\ref{lemmaconj3}.] First note that  $S_i+N/N\ (i=1, 2)$ are two Levi subalgebras of the Lie algebra $L=M/N$. By the Levi--Malcev theorem for~$M/N$, there is a derivation $\delta\in {\rm Der}(L)$ such that $\exp(\delta)$ maps $S_1+N/N$ onto $S_2+N/N$. Because $[NM]=0, \delta$ lifts to a derivation of~$M$. So in proving the theorem, we may, and shall, assume that $S_1+N=S_2+N$.

Arguing in the same way as the discussion preceding Example~\ref{example3.1}, we know that
$S_2 = \{(x, c(x))\, | \, x \in S_1\}$ for some $c\in {\rm Hom}_{S_1}(S_1, N)$. We will show that the linear map $f: M\rightarrow M$ def\/ined by
\begin{gather*}
f(x)=\begin{cases} c(x), & x\in S_1, \\
  0  ,&  x\in B,
\end{cases}
\end{gather*}
has the desired properties. $f$ is well-def\/ined  because $M=B\oplus S_1$, and since $f(M)\subseteq N\subseteq B$ then $f^2=0$. To see that $f$ is a derivation, note that for $x,y\in S_1$, $a,b\in B$, we have
\begin{gather*}
f([(x+a)(y+b)]) = f([xy])  =  c([xy]), \\
[(x+a)f(y+b){]}+[(f(x+a))(y+b)]  =  [(x+a)c(y)]= c([xy]),
\end{gather*}
where we used $[ac(y)] \in [BN]=0$ to ensure the last equality.

Since $f^2=0$, we have $\exp(f)=Id+f$, whence $\exp(f)(x)=x+c(x)$ for $x \in S_1$.
Consequently, $\exp(f)$ maps $S_1$ onto $S_2$. This completes the proof of the theorem.
\end{proof}

\section{Left central Leibniz algebras} \label{section4} 
In this section, $M$ is a left central Leibniz algebra with Leibniz kernel $N=C(M)$.  Recall that this means that $N \subseteq Z(M)$, i.e.~$[MN]=0$.
Note that $M$ satisf\/ies the hypotheses of both Theorem~\ref{lemmaconj1} and Theorem~\ref{lemmaconj3}, so that the Levi--Malcev theorem holds in~$M$.  This suggests that left central Leibniz algebras are more manageable than  general Leibniz algebras.

We introduce and study a very useful bilinear pairing attached to~$M$.
For $a, b \in M$, write
\begin{gather*}
\psi_M(a, b) = \psi(a, b) := [a b]+[ba].
\end{gather*}
Then $\psi: M \times M \rightarrow N$ is a symmetric bilinear map.
Even though $\psi$ takes values in $N$, it still makes sense to use terminology associated with the more familiar situation of a $\mathbb{C}$-valued bilinear form on $M$.
Thus $a \in M$ is \emph{isotropic} if $\psi(a, a)=0$ ($[aa]=0$); otherwise,
$a$ is \emph{anisotropic}. $a$ and $b$ are
\emph{orthogonal} if $\psi(a, b)=0$. For a subspace $U\subseteq M$, def\/ine
\begin{gather*}
  U^{\perp} := \{ b \in M \, | \, \psi(a, b)=0, \  a \in U \}, \\
{\rm rad}(U) := U \cap U^{\perp}.
\end{gather*}
$U$ is \emph{totally isotropic} if $U\subseteq U^{\perp}$ (this is equivalent to every element of~$U$ being isotropic, and also to $U= {\rm rad}\, U$), and  \emph{nondegenerate} if ${\rm rad}\, U=0$. $U, V \subseteq M$ are \emph{orthogonal} if $V\subseteq U^{\perp}$.
The \emph{radical} of~$\psi$ (or~$M$) is
\begin{gather*}
R := {\rm rad}(M) = M^{\perp}= \{x \in M \, | \, \psi(x, a)=0, \ a \in M\}.
\end{gather*}

These def\/initions apply to all Leibniz algebras. $\psi$ is particularly useful in the study of central Leibniz algebras because of the next result.
\begin{lemma} Suppose that~$M$ is a left central Leibniz algebra. Then $\psi$ is  associative  in the sense that for $a, b, c \in M$ we have
\begin{gather}\label{Rdef}
\psi([ab], c) = \psi(a, [bc]).
\end{gather}
In particular, $R$ is a $2$-sided ideal of $M$.
\end{lemma}

\begin{proof} The identities $([ab]+[ba])c = a([bc]+[cb])=0$ hold for all
$a, b, c \in M$. Therefore,
\begin{gather*}
\psi([ab], c)  = [[ab]c]+[c[ab]] = -([[ba]c]+[c[ba]])
 = -([b[ac]] -[a[bc]]+[[cb]a]+[b[ca]]) \\
\hphantom{\psi([ab], c)}{}
= [a[bc]]+[[bc]a] = \psi(a, [bc]).
\end{gather*}
This proves (\ref{Rdef}). That $R$ is a $2$-sided ideal is an immediate consequence.
\end{proof}

  Because $M$ is a left central Leibniz algebra, we have $N \subseteq R$. Thus,
  $M/R$ is a Lie algebra equipped with the symmetric nondegenerate pairing induced by~$\psi$.
  Introduce the poset (with respect to containment)
  \begin{gather*}
\mathcal{L} := \{\mbox{Lie subalgebras of $M$}\}.
\end{gather*}

  \begin{lemma}\label{lemL+R} Let $L \subseteq M$ be a Leibniz subalgebra. Then the following are equivalent:
  \begin{gather*}
 (a)\ L \in \mathcal{L}; \\
 (b)\ L \ \mbox{is a totally isotropic subspace}; \\
 (c)\ L+R \in \mathcal{L}.
\end{gather*}
  \end{lemma}

\begin{proof} The equivalence of (a) and (b) follows directly from the def\/initions. Next, if $L+R$ is a Lie algebra then $L$ is necessarily a Lie subalgebra, so that (c) $\Rightarrow$ (a). Conversely, if~(a) holds then~$L$ is a Lie subalgebra of~$M$. In this case, $L+R$ is itself a Leibniz subalgebra because $R$ is a $2$-sided ideal.
 Moreover, if $a\in R$, $b\in L$, then
\begin{gather*}
\psi(a+b, a+b) = \psi(a, a)+2[bb]+2\psi(a, b) = 0.
\end{gather*}
Therefore, $R+L$ is totally isotropic, and the equivalence of (a) and (b) applied to
$R+L$ shows that (c) holds. This completes the proof of the lemma.
\end{proof}

  \begin{lemma}\label{lemisoideal} Suppose that $U \subseteq M$ is a totally isotropic $2$-sided
  ideal of~$M$.
  Then $U \in \mathcal{L}$ and $U' \subseteq R$.
  \end{lemma}

 \begin{proof}
 Because $U$ is a totally isotropic Leibniz subalgebra of $M$, Lemma~\ref{lemL+R} tells us that $U \in \mathcal{L}$. Now let $a, b \in U$, $x \in M$. Then $[b x] \in U$, and the associativity of $\psi$ (Lemma~\ref{Rdef}) together with the fact that~$U$ is totally isotropic implies that $\psi([ab], x) = \psi(a, [bx])=0$.  Therefore $[ab] \in R$ for all $a, b \in U$, i.e.\
$U'\subseteq R$.
This completes the proof of the lemma.
\end{proof}

\begin{lemma}\label{lemL/R} Suppose that $R\subseteq L \subseteq M$ such that $L/R \subseteq M/R$ is a~semisimple  Lie subalgebra.
Then $L \in \mathcal{L}$.
\end{lemma}

\begin{proof}
$L$ is clearly a Leibniz subalgebra of $M$. Let $S \subseteq L$ be a
Levi subalgebra. Because $L/R$ is semisimple we have $L=R+S$, and this is
a Lie subalgebra thanks to the equivalence of~(a) and~(c) in Lemma~\ref{lemL+R}.
\end{proof}

\begin{lemma}\label{solvideal} $M/R$ has  no  nonzero semisimple ideals.
In particular, every minimal ideal of~$M/R$ is abelian.
\end{lemma}

\begin{proof}
Suppose that $U/R$ is an ideal of $M/R$ that is also a semisimple
Lie algebra. By Lemma~\ref{lemL/R}, $U\in \mathcal{L}$, and in particular it is totally isotropic.
Then
$U' \subseteq R$ by Lemma~\ref{lemisoideal}. Therefore~$U/R$ is both abelian and semisimple, whence it reduces to~$0$.
This proves the f\/irst assertion of the lemma. Because every minimal  ideal
of the Lie algebra~$M/R$ is either abelian or semisimple, the second assertion follows. This completes the proof of the lemma.
\end{proof}

\begin{proposition}\label{prop4.1} Let $\mathcal{L}^*$ be the set of  maximal elements  of~$\mathcal{L}$. Then
\begin{gather*}
R = \bigcap_{L \in \mathcal{L}^*} L.
\end{gather*}
\end{proposition}

\begin{proof}
Let $L \in \mathcal{L}^*$. Then $R+L$ is also a Lie subalgebra
by Lemma~\ref{lemL+R}, so that $R\subseteq L$ because $L$ is maximal.
  This shows that $R$ is contained in every $L \in \mathcal{L}^*$, and hence also in their intersection.

To prove the opposite containment, we use induction on $\dim N$. There is nothing to prove if~$M$ is a Lie algebra, so we may, and shall, assume that this is not the case. Therefore, $N \not= 0$. Suppose that $\dim N \geq 2$.
Then every hyperplane $N_0 \subseteq N$ is a $2$-sided ideal in $M$, and is itself contained in every maximal Lie subalgebra $L$. Then by induction we obtain
\begin{gather*}
\bigcap_{L \in \mathcal{L}^* } L/N_0 = R(M/N_0),
\end{gather*}
which says that if $a \in \cap L$ then $\psi(a, x) \in N_0$ for all $x \in M$. Then
$\psi(a, x) \in \cap N_0 =0$ (the last intersection ranging over hyperplanes of $N$),
and thus $a \in R$, as required
This reduces us to the case when $\dim N=1$, so that we can think of $\psi$ as a
$\mathbb{C}$-valued bilinear form. We assume this for the remainder of the proof.

Because $M$ is not a Lie algebra, $R \not= M$. Suppose that
$R$ has codimension $1$. Then every nonzero element of $M/R$ is anisotropic,
so that if $a \in M\setminus{R}$ then  $a$ cannot be contained
in a~Lie subalgebra of $M$. Thus in this case, $R$ is  the \emph{unique} element in $\mathcal{L}^*$,
and the desired result is clear.

Finally, suppose that $M/R$ has dimension at least~$2$. Because~$\psi$ is a nondegenerate
$\mathbb{C}$-valued bilinear form on~$M/R$, $M/R$ is \emph{spanned}
by \emph{isotropic} vectors, i.e., elements $a+R$ with $[aa]=0$.  Such elements $a$ are contained in some maximal Lie subalgebra; therefore, if $b \in \cap_{L \in \mathcal{L}^*} L$,
$a$~and~$b$ generate a Lie subalgebra of~$M$, so that
$[ab]+[ba]=0$. Because the isotropic vectors $a+R$ span, we can conclude that
$[ab]+[ba]=0$ for all $a \in M$, whence $b \in R$. This completes the proof of the proposition.
\end{proof}

We will need the next lemma in Section~\ref{section5}.
\begin{lemma}\label{lemBperp} Let $R\subseteq B\subseteq M$ satisfy $B/R:= B(M/R)$.
Then $B^{\perp} \subseteq B$.
\end{lemma}

\begin{proof} Write
$B^{\perp} = S_0\oplus B_0$, where $S_0$ and $B_0$
 are a Levi factor and the solvable radical respectively for~$B^{\perp}$.
Because~$B$ is an ideal in~$M$ then so too is~$B_0$. Therefore,
$B_0 \subseteq B\cap B^{\perp}$. Thus $B_0 = {\rm rad}(B)$, and in particular $B_0$ is totally isotropic.
So we see that there is an orthogonal decomposition $B^{\perp}=S_0 \perp B_0$,
and since both summands are totally isotropic
then~$B^{\perp}$ is a totally isotropic ideal of~$M$.
Now we can apply Lemma~\ref{lemisoideal}, with $B^{\perp}$ playing the r\^{o}le of $U$, to conclude that $S_0 \subseteq R$. Therefore, $B^{\perp} \subseteq B$, and the lemma is proved.
\end{proof}

We now describe the construction of a class of left central Leibniz algebras that depends
only on Lie-theoretic data $(L_1, L_2, R, \alpha, \pi)$ satisfying (a)--(c) below. The set-up is as follows:

(a)
a pair of Lie algebras $L_1$, $L_2$ with a common ideal $R$ and $L_2/R$ abelian:
\begin{gather*}
 0 \rightarrow R \rightarrow L_1  \rightarrow L_1/R \rightarrow 0 \\
  \ \ \ \ \ \ \, \| \\
 0 \rightarrow R \rightarrow L_2 \rightarrow L_2/R \rightarrow 0
\end{gather*}
Set $Z:= R\cap Z(L_1)\cap Z(L_2)$.

(b)\ a morphism of Lie algebras $ \alpha: L_1 \rightarrow {\rm Der}(L_2)$
with $\alpha|_R = {\rm ad}_R$, $\alpha(L_1)|_R={\rm ad}_{L_1}|_R$. Setting
$\alpha(x_1)(y_1)=x_1.y_1$, these assumptions amount to
\begin{gather*}
[x_1x_2].y_1  =  x_1.(x_2.y_1)-x_2.(x_1.y_1), \\
x_1.[y_1y_2]  =  [(x_1.y_1)y_2]+[y_1(x_1.y_2)], \\
x_1.y_1  =  [x_1y_1] \qquad (x_1\in R\  \mbox{or} \ y_1 \in R),
\end{gather*}
for $x_1, x_2 \in L_1$, $y_1, y_2 \in L_2$, and where $[ \ ]$ is the bracket in $L_1$ or $L_2$.

(c)\ an \emph{injective} morphism of left $L_1$-modules $\pi: L_1/R \rightarrow {\rm Hom}_{\mathbb{C}}(L_2/R, Z)$ such that ${\rm im}\, \pi$ annihilates no nonzero elements of~$L_2/R$. ($L_1$ acts on~$L_2/R$ via $\alpha$, trivially on~$Z$, and with the induced left action
on ${\rm Hom}_{\mathbb{C}}(L_2, Z)$). Lifting~$\pi$ to a morphism of left $L_1$-modules
$\pi: L_1\rightarrow  {\rm Hom}_{\mathbb{C}}(L_2, Z)$ and setting $\psi': L_1\times L_2 \rightarrow Z, (x_1, y_1)
\mapsto \pi(x_1)(y_1)$, these assumptions mean that
\begin{gather*}
  \psi'([x_1x_2], y) = \psi'(x_1, x_2.y), \\
  R = \{ x_1 \in L_1\, | \, \psi'(x_1, L_2)=0\}=\{y_1 \in L_2\, | \, \psi'(L_1, y_1)=0\}.
\end{gather*}

Def\/ine $M = (L_1\oplus L_2, [ \ ])$, where
\begin{gather}\label{Mconstruct}
[(x_1,  y_1)(x_2, y_2)] := ([x_1x_2], x_1.y_2-x_2.y_1+[y_1y_2]+\psi'(x_2, y_1)).
\end{gather}
Thus $L_1$, $L_2$ are naturally Lie subalgebras of $M$. One calculates that $M$ is a left Leibniz algebra, $L_2$~is a $2$-sided ideal, and $Z\oplus Z\subseteq Z(M)$. Moreover,
$\psi((x_1, y_1), (x_2, y_2)) = (0, \psi'(x_2, y_1)+\psi'(x_1, y_2)) \in Z\oplus Z$.
So $M$ is a left central Leibniz algebra with radical $R\oplus R$.
Furthermore, $D :=\{ (a, -a) \, | \, a \in R\}$
 is a $2$-sided ideal. So $\tilde{M} := M/D$ is itself a central Leibniz algebra,
 and the radical $R\oplus R/D$ of the induced bilinear form $\tilde{\psi}$ on $\tilde{M}$ is
 naturally identif\/ied with $R$. We state these conclusions as
 \begin{lemma}\label{lemtildeM}  Suppose that $(L_1, L_2, R, \alpha, \pi)$ satisfies
 assumptions $(a)$--$(c)$. Then
 $\tilde{M}= \tilde{M}(L_1, L_2,$ $R, \alpha, \pi)$ is a left central Leibniz algebra with radical $R$.
 \end{lemma}

\section{Left central Leibniz algebras of rank 1}\label{section5}

For a Leibniz algebra $M$, we def\/ine the \emph{rank} of $M$
to be the \emph{dimension} of $N=C(M)$, and denote it by ${\rm rk}(M)$.  This is a  useful
invariant for left central Leibniz algebras, because in this case any subspace $N_0 \subseteq N$ is a
$2$-sided ideal of~$M$ and~$M/N_0$ is a left central Leibniz algebra with ${\rm rk}(M/N_0)=\dim N - \dim N_0$.  In this way, we can try to reduce questions about left central Leibniz algebras to those
of rank~$1$.  These are easier to deal with, because the form $\psi$ may be considered to be a \emph{trace form} in the usual sense, i.e.\ a $\mathbb{C}$-valued associative bilinear form.  This method was already used in the proof of Proposition~\ref{prop4.1}.

If $M_1$, $M_2$ are two Leibniz algebras, their orthogonal sum
$M_1 \perp M_2$ is the direct sum $M_1\oplus M_2$ with product $[(a_1, b_1)(a_2, b_2)]=([a_1a_2], [b_1b_2])$. Then $M_1\perp M_2$ is a Leibniz
algebra, $M_1$ and $M_2$ are orthogonal ideals (with respect to $\psi_{M_1\perp M_2}$),
and ${\rm rk}(M_1\perp M_2)={\rm rk}(M_1)+{\rm rk}(M_2)$. Thus $rk$ is additive over orthogonal sums.

If $M$ has index $r \geq 1$, we can f\/ind $r$ hyperplanes $N_1, \dots, N_r$ of $N$ such that $\cap_{i=1}^r N_i = 0$. Then there is an injective morphism of Leibniz algebras
\begin{gather*}
M \rightarrow \perp_{i=1}^r M/N_i, \qquad a \mapsto (a+N_1, \dots, a+N_r).
\end{gather*}
In this way, any left central Leibniz algebra of rank $r\geq 1$ is isomorphic to a subalgebra of an
orthogonal  sum of $r$ left central Leibniz algebras of rank $1$.

\begin{theorem}\label{thm5.1} Suppose that $M$ is a  left central Leibniz algebra of rank~$1$, and let
$B$ be as in Lemma~{\rm \ref{lemBperp}}.
Then there is
 at least one  maximal Lie algebra $L \in \mathcal{L}^*$ with the following properties:
\begin{gather*}
 (a)\ \mbox{$L$ is a  maximal  isotropic subspace of $M$}; \\
 (b)\ \mbox{$L/R$ is nilpotent};\\
 (c)\ L \subseteq B.
\end{gather*}
\end{theorem}

\begin{proof} If $L \in \mathcal{L}^*$ then $R \subseteq L$ (cf.\ Proposition~\ref{prop4.1}).
 So part~(b) makes sense.
To prove the theorem, we use induction on $\dim M$. If every nonzero element of~$M/R$ is
anisotropic, then~$R$ is both the unique maximal Lie subalgebra of $M$ and the unique maximal totally isotropic subspace, in which case all parts of the theorem are obvious.
So we may, and shall, assume that~$M/R$ contains nonzero isotropic elements.  As in the proof of Proposition~\ref{prop4.1}, this implies that $\dim (M/R)\geq 2$.

We assume f\/irst that there is a (nonzero) ideal $U/R$ in $M/R$ which is
totally isotropic. If $M/R$ is solvable, we choose any such ideal $U/R$;
 if $M/R$ is nonsolvable we take $U = B^{\perp}$. Note that in the second case,
 $B$ is a proper ideal of~$M$, so that $B^{\perp}/R$ is nonzero thanks to the assumption
 that~$M$ has rank~$1$. Moreover, $B^{\perp}/R$ is totally isotropic in the second case because of Lemma~\ref{lemBperp}.

With $U/R$ chosen in this way,  $U^{\perp}$
is a proper ideal of $M, U^{\perp} \subseteq B$, and $U=$ rad$(U^{\perp})$. By induction there is
a Lie subalgebra $L \subseteq U^{\perp}$ which is a maximal isotropic subspace of $U^{\perp}$
 and such that $L/U$ is nilpotent.  If $L_1$ is a maximal isotropic subspace of $M$ that contains~$L$, then $L_1 \subseteq U^{\perp}$, and this implies that $L=L_1$ since $L$ is maximal isotropic in $U^{\perp}$. So $L$ is maximal isotropic in~$M$.

It remains to show that $L/R$ is nilpotent. Let $a \in M, u \in U, x \in L$. Then
$\psi(a, [ux]) = \psi([au], x)=0$, which shows that $[LU] \subseteq R$. Another way to say this is
$U/R \subseteq Z(L/R)$. Since we already know that $L/U$ is nilpotent, we can conclude that
$L/R$ is too. This completes the proof of the theorem in this case.
 Thus we may, and now shall, suppose that \emph{no}
nonzero  ideal~$U/R$ of~$M/R$ is  totally isotropic.

Now take any minimal nonzero ideal $U/R\subseteq M/R$. Then $U/R$ is nondegenerate, and because~$M$ has rank $1$ there is an orthogonal
decomposition $M/R = U/R \perp U^{\perp}/R$. Continuing in this way, we obtain
an orthogonal decomposition
\begin{gather*}
M/R = U_1/R \perp \cdots \perp U_k/R
\end{gather*}
with each $U_j/R$ a minimal, nondegenerate ideal of $M/R$.
By Lemma~\ref{solvideal}, each $U_j/R$ is \emph{abelian}. Then $M/R$ is itself an abelian Lie algebra, and every nonzero element of $M/R$ generates a~$1$-dimensional ideal. In particular,
 since $M/R$ contains nonzero isotropic elements, it also has a~nonzero  isotropic  ideal. This contradiction completes the proof of the theorem.
 \end{proof}

\section{Symmetric Leibniz algebras}\label{section6}

\begin{lemma}\label{lemLsymm} Suppose that $M$ is a left Leibniz algebra. Then $M$ is a symmetric Leibniz algebra if, and only if, $M'\subseteq R$.
\end{lemma}

\begin{proof}
 In any left Leibniz algebra we have
\begin{gather*}
 [a[bc]]+[[ac]b] = [[ab]c]+[b[ac]]+[[ac]b] = [[ab]c]+\psi([ac], b).
 \end{gather*}
Therefore, $M$ is a symmetric Leibniz algebra if, and only if,
$\psi([ab], c)=0$ for all $a, b, c \in M$. The assertion of the lemma is just a restatement of this.
\end{proof}

Consider the left central Leibniz algebra $M=M(L_1, L_2, R, \alpha, \pi)$ with product~(\ref{Mconstruct}) and its quotient
algebra $\tilde{M}=M/D$ (cf.\ Lemma~\ref{lemtildeM}) introduced at the end of Section~\ref{section4}. Since $D$ is contained in the radical of $M$, it follows from Lemma~\ref{lemLsymm} that $\tilde{M}$ is symmetric if, and only if, $M$ is symmetric. From~(\ref{Mconstruct}) and Lemma~\ref{lemLsymm} once more, this holds if, and only if, $L_1/R$ is abelian (recall that $L_2/R$ is abelian by construction) and $x_1.y_2-x_2.y_1 \in R$ $(x_1, x_2 \in L_1, y_1, y_2 \in L_2)$. We assert that the second condition is a consequence of the f\/irst. Indeed, because $\pi: L_1/R \rightarrow {\rm Hom}_{\mathbb{C}}(L_2/R, Z)$ is an injection of $L_1$-modules then $L_1$ acts trivially on ${\rm im}\,\pi$, i.e., ${\rm im}\,\pi$ annihilates each element $x_1.y_1$.
 Since ${\rm im}\,\pi$ has no nonzero invariants in its action on~$L_2/R$, it follows that $x_1.y_1 \in R$, and since this holds for any $x_1\in L_1$, $y_1\in L_2$ then the second condition indeed follows, as asserted. Thus we have proved

 \begin{lemma} The left central Leibniz algebras $M(L_1, L_2, R, \alpha, \pi)$ and
 $\tilde{M}=M/D$ are symmetric Leibniz algebras if, and only if, $L_1/R$ and $L_2/R$ are both abelian Lie algebras.
 \end{lemma}

Now suppose that $M_0$ is a symmetric Leibniz algebra of rank $1$.
Thus $N= C(M_0)$ has dimension~$1$ and~$\psi$ def\/ines a
trace form on~$M_0$ with radical~$R$. Consequently, we can f\/ind a pair of
maximal, totally isotropic subspaces $L_1, L_2 \subseteq M_0$
such that $L_1\cap L_2 = R$ and $L_1+ L_2$ has codimension at most $1$ in~$M_0$. \
($L_1+L_2=M_0$ if, and only if, $\dim(M_0/R)$ is \emph{even}.)\
 Because $M_0/R$ is an abelian Lie algebra by Lemma~\ref{lemLsymm}, it follows that
 $H :=L_1+L_2$ is a $2$-sided ideal with Lie subalgebras~$L_1$,~$L_2$;
 indeed, $L_1, L_2 \in \mathcal{L}^*$.

We will show that $H \cong \tilde{M}(L_1, L_2, R, \alpha, \pi)$ for suitably def\/ined maps $\alpha$, $\pi$. Both $L_1$, $L_2$ are ideals of~$M_0$, in particular the left adjoint def\/ines
a morphism of Lie algebras $\alpha: L_1 \rightarrow {\rm Der}(L_2)$ satisfying the assumptions
of~(b) at the end of Section~\ref{section4}. Similarly, $R= {\rm  rad} (M)$ and we def\/ine the morphism of
$L_1$-modules $\pi: L_1/R\rightarrow {\rm Hom}_{\mathbb{C}}(L_2/R, N)$ as $\pi(x_1, y_1):= \psi(x_1, y_1)$ $(x_1 \in L_1$, $y_1 \in L_2)$. It is easily seen that the assumptions of~(c) at the end of Section~\ref{section4} also hold. (The main point is  the associativity of~$\psi$.)  Thus the set-up discussed in Section~\ref{section4} holds, and we can apply Lemma~\ref{lemtildeM} to obtain the left central Leibniz algebras
$M(L_1, L_2, R, \alpha, \pi)$ and~$\tilde{M}$. We assert that
$H \cong \tilde{M}$. To see this, one checks that the map
\begin{gather*}
\varphi: \ M(L_1, L_2, R, \alpha, \pi) \rightarrow H, \qquad (x_1, y_1) \mapsto x_1+y_1,
\end{gather*}
 is a surjective morphism of Leibniz algebras with kernel
$\{(x_1, y_1)\, | \, x_1+y_1 = 0\} = \{(a, -a)\, | \, a \in R\}$.
Now our assertion follows from the very construction of~$\tilde{M}$.
We have proved

\begin{theorem}\label{thm6.3} A symmetric Leibniz algebra of rank~$1$
has an ideal of codimension at most~$1$ isomorphic to
$\tilde{M}(L_1, L_2, R, \alpha, \pi)$.
\end{theorem}

\subsection*{Acknowledgements}
This work was supported by the NSF (G.M.) and the
Simons Foundation (G.Y.). The authors thank the (anonymous) referees for helpful comments.

\pdfbookmark[1]{References}{ref}
\LastPageEnding

\end{document}